\documentclass{amsart}

\usepackage{amssymb}
\usepackage{calrsfs}
\usepackage{enumerate}

\newtheorem{theorem}{Theorem}[section]
\newtheorem{cor}[theorem]{Corollary}
\newtheorem{lemma}[theorem]{Lemma}

\theoremstyle{definition}

\theoremstyle{remark}

\numberwithin{equation}{section}

\def\cM{{\mathcal M}}
\def\cN{{\mathcal N}}

\def\cA{{\mathcal A}}

\newcommand{\vp}{\varphi}

\def\cB{{\mathcal B}}

\begin{document}

\title{Zero product preservers of C*-algebras}

\author{Ngai-Ching Wong}

\date{January 1, 2007; to appear in the ``Proceedings of the Fifth Conference on Function Space'',
Contemporary Math.}



\dedicatory{Dedicated to Professor Bingren Li on the occasion of his 65th birthday (1941.10.7 -- )}

\address{Department of Applied Mathematics, National Sun Yat-sen
University,  and National Center for Theoretical Sciences,
  Kaohsiung, 80424, Taiwan, R.O.C.} \email{wong@math.nsysu.edu.tw}

\thanks{This work is partially supported by Taiwan NSC grant 95-2115-M-110-001-}

\begin{abstract}
Let $\theta: \cA\to \cB$ be a zero-product preserving bounded linear map between
C*-algebras. Here neither $\cA$ nor $\cB$ is necessarily unital. In this note, we investigate
when $\theta$ gives rise to a Jordan
homomorphism.  In particular, we
show that $\cA$ and $\cB$ are
isomorphic as Jordan algebras if $\theta$ is bijective and sends
zero products of self-adjoint elements to zero products.
They are isomorphic as C*-algebras if  $\theta$ is bijective and preserves the full zero product structure.
\end{abstract}

\maketitle

\section{Introduction}
Let $\cM$ and $\cN$ be algebras over a field $\mathbb F$ and $\theta
: \cM\to \cN$ a linear map. We say that $\theta$ is a
\emph{zero-product preserving} map if $\theta(a)\theta(b)=0$ in
$\cN$ whenever $ab=0$ in $\cM$.  The canonical form of a linear zero product preserver,
 $\theta=h\vp$,  arises
from  an element $h$
in the center of $\cN$ and an algebra
homomorphism $\vp: \cM\to \cN$. In
\cite{CKLW}, we see that in many interesting cases
zero-product preserving linear maps arise in this way.

We are now interested in the $C^*$-algebra case.
There are 4 different versions of zero products:
$ab=0$, $ab^*=0$, $a^*b=0$ and $ab^*=a^*b=0$.
Surprisingly, the original version $ab=0$ is the least, if any, geometrically meaningful,
while the others mean $a,b$ have orthogonal initial spaces, or orthogonal range spaces, or both.
Using the orthogonality conditions,
 the author showed in \cite{wong-seam} that
a bounded linear  map $\theta:\cA\to\cB$ between
C*-algebras is a triple homomorphism if and only if $\theta$ preserves the fourth disjointness
 $ab^*=a^*b=0$ and $\theta^{**}(1)$ is a partial isometry.  Here,
 the triple product of a C*-algebra is defined by $\{a,b,c\}=(ab^*c + cb^*a)/2$,
 and $\theta^{**}:\cA^{**}\to\cB^{**}$ is the
 bidual map of $\theta$.
 See also \cite{ABEV} for a similar result dealing with the case $ab=ba=0$.
We shall deal with the first and original case in this note.  The other cases will be dealt with in a subsequent paper.

There is a common starting point of all these 4 versions.  Namely, we can consider first the zero products
$ab=0$ of self-adjoint elements $a,b$ in $\cA_{sa}$.
In \cite{wolff} (see also \cite{schweizer}), Wolff shows that if
$\theta:\cA\to \cB$ is a bounded linear map  between unital
$C^*$-algebras preserving the involution and zero products of
self-adjoint elements in $\cA$ then $\theta=\theta(1)J$ for a Jordan
$*$-homomorphism $J$ from $\cA$ into $\cB^{**}$.
In \cite{CKLW}, the involution preserving assumption
is successfully removed.
Modifying the arguments in \cite{CKLW}, we will further relax
the condition that the C*-algebras are unital in this note.
In particular, we
show that $\cA$ and $\cB$ are
isomorphic as Jordan algebras if  $\theta$ is bijective and sends self-adjoint elements with zero products in $\cA$ to
elements (not necessarily self-adjoint, though) with zero products in $\cB$.
They are isomorphic as C*-algebras if $\theta$ is bijective and preserves the full zero product structure.

\section{Results}

In the following, $\cA,\cB$ are always C*-algebras not necessarily with identities.
$\cA_{sa}$ denotes the (real) Jordan-Banach algebra consisting of all self-adjoint elements of $\cA$.

Recall that a linear map $J$ between two algebras is said to be a
\emph{Jordan homomorphism} if $J(xy +yx) = J(x)J(y)+J(y)J(x)$ for
all $x,y$. If the underlying field has characteristic
not $2$, this condition is equivalent to that $J(x^2)=(Jx)^2$ for
all $x$ in the domain.  We also have the identity $J(xyx)=J(x)J(y)J(x)$ for all $x,y$ in this case.

\begin{lemma}\label{lem:new1}
Let $J : \cA_{sa} \longrightarrow \cB$ be a bounded Jordan homomorphism.  Then $J$
 sends zero products in $\cA_{sa}$ to zero products in
$\cB$.
\end{lemma}
\begin{proof}
Let $a,b$ be self-adjoint elements in $\cA$ and $ab=0$.
We want to prove that $J(a)J(b)=0$.  Without loss of generality, we can
assume that  $a\geq 0$.
Let $a'$ in $A_{sa}$ satisfy that ${a'}^2=a$.  We have
$a'b=0$.  By the identities
$
0 = J(a'ba')
 =J(a')J(b)J(a')
$
and
$
0 =J(a'b+ba')=J(a')J(b)+J(b)J(a')
$,
we have
$
J(a)J(b)=J({a'}^2)J(b)=J(a')^2J(b)=0.
$
\end{proof}

Recall that when we consider
$\cA^{**}$ as the enveloping W*-algebra of $\cA$,
the multiplier algebra $M(\cA)$ of $\cA$ is the C*-subalgebra of $\cA^{**}$,
$$
M(\cA) =\{x\in \cA^{**} : x\cA \subseteq \cA \text{ and } \cA x\subseteq \cA\}.
$$
Elements in $M(\cA)_{sa}$ can be approximated by both monotone increasing and decreasing bounded nets from
$\tilde\cA_{sa} = \cA_{sa} \oplus \mathbb{R} 1$ (see, e.g., \cite{brown1988}).
In case $\cA$ is unital, $M(\cA)=\cA$.

\begin{lemma}\label{lem:ma}
Let $\theta : \cA_{sa}\to \cB$ be a bounded linear map sending zero products in $\cA_{sa}$ to zero products in $\cB$.
Then the restriction of $\theta^{**}$ induces a bounded linear map, denoted again by $\theta$, from $M(\cA)_{sa}$ into $\cB^{**}$, which sends
zero products in $M(\cA)_{sa}$ to zero products in $\cB^{**}$.
\end{lemma}
\begin{proof}
First we consider the case $b\in \cA_{sa}$, and $p$ is an open projection in $\cA^{**}$ such that $pb=0$.
For any self-adjoint element $c$ in the hereditary C*-subalgebra $h(p)=p\cA^{**}p\cap \cA$ of $\cA$, we have
$cb=0$ and thus $\theta(c)\theta(b)=0$.  By the weak* continuity of $\theta^{**}$, we have
$\theta^{**}(p\cA^{**}_{sa}p)\theta(b)=0$.  In particular, $\theta^{**}(p)\theta(b)=0$.

Let $a,b$ be self-adjoint elements in $M(\cA)$ with $ab=0$.
We want to prove that $\theta(a)\theta(b)=0$.
Without loss of generality, we can assume both $a,b$ are positive.
Let $0\leq a_\alpha +\lambda_\alpha \uparrow a$ be a monotone increasing net from $\tilde\cA_{sa}$.
Since $0\leq b(a_\alpha +\lambda_\alpha) b \uparrow bab=0$,
we have $(a_\alpha +\lambda_\alpha) b =0$ for all $\alpha$.
Similarly, there is a monotone increasing net $0\leq b_\beta +s_\beta\uparrow b$ from $\tilde\cA_{sa}$ such that
$(a_\alpha +\lambda_\alpha)(b_\beta +s_\beta)=0$ for all $\beta$.  We can assume the real scalar
$\lambda_\alpha \neq 0$.  Then $s_\beta=0$ for all $\beta$.
In particular, we see that $a_\alpha$ commutes with all $b_\beta$.
In the abelian C*-subalgebra of $M(\cA)$ generated by $a_\alpha$, $b_\beta$ and $1$,
we see that $a_\alpha  +\lambda_\alpha$ can be approximated in norm by finite real linear combinations of open projections
disjoint from $b_\beta$.  By the first paragraph, we have
$\theta(a_\alpha +\lambda_\alpha)\theta(b_\beta)=0$.

 By the weak* continuity of $\theta^{**}$ again, we see
that  $\theta(a_\alpha + \lambda_\alpha)\theta(b) =
 \lim_\beta \theta(a_\alpha +\lambda_\alpha)\theta(b_\beta) = 0$ for each $\alpha$, and then
$\theta(a)\theta(b) = \lim_\alpha \theta(a_\alpha +\lambda_\alpha)\theta(b) = 0$.
\end{proof}

With Lemma \ref{lem:ma}, results in \cite{CKLW} concerning zero product preservers of unital C*-algebras can be extended easily to
the non-unital case.
We  restate \cite[Lemmas 4.4 and 4.5]{CKLW} below, but now here $\cA$ does not necessarily have an identity.

\begin{lemma}\label{lem:commutant}
Let $\theta: \cA\to \cB$ be a bounded linear map sending zero products in $\cA_{sa}$ to zero products in $\cB$.
For any $a$ in $M(\cA)$, we have
\begin{enumerate}[\rm(i)]
    \item $\theta(1)\theta(a)=\theta(a)\theta(1)$,
    \item $\theta(1)\theta(a^2)=(\theta(a))^2$.
\end{enumerate}
If $\theta(1)$ is invertible then $\theta = \theta(1)J$ for a bounded Jordan
homomorphism $J$ from $\cA$ into $\cB$.
\end{lemma}

\begin{theorem}\label{theorem:4.6}
Two C*-algebras $\cA$ and $\cB$
are isomorphic as Jordan algebras if and only
if there is a bounded bijective linear map $\theta$ between them sending
zero products in $\cA_{sa}$ to zero products in
$\cB$.
If $\theta$ is just surjective, then $\cB$ is isomorphic to the  C*-algebra $\cA/\ker \theta$ as Jordan algebras.
\end{theorem}
\begin{proof}
One way follows from Lemma \ref{lem:new1}.
Conversely, suppose $\theta(\cA)=\cB$.
Since $\theta(1)\theta(a^2)=\theta(a)^2$ for all $a$ in $\cA$ and $\cB=\cB^2$,
we have $\theta(1)\cB=\cB$.  Thus,
the central element $\theta(1)$ is invertible.  Lemma \ref{lem:commutant} applies,
by noting that closed Jordan ideals of C*-algebras are two-sided ideals \cite{Civin65}.
\end{proof}

In case $\theta$ preserves all zero products in $\cA$, we have the
following non-unital version of \cite[Theorem 4.11]{CKLW}.

\begin{theorem}\label{theorem:calgfull}
Let $\theta$ be a surjective bounded linear map from a
$C^*$-algebra $\cA$ onto a $C^*$-algebra $\cB$. Suppose that
$\theta(a)\theta(b)=0$ for all $a,b\in\cA$ with $ab=0$. Then
$\theta(1)$ is a central element and invertible in
$M(\cB)$. Moreover, $\theta=\theta(1)\vp$ for a surjective algebra
homomorphism $\vp$ from $\cA$ onto $\cB$.
\end{theorem}
\begin{proof}
First, we have already seen in the proof of Theorem \ref{theorem:4.6} that $\theta(1)$ is a central
element and invertible in $M(\cB)$.  Second, we observe that to utilize the results
 \cite[Theorems 4.12 and 4.13]{CKLW} of Bre\v{s}ar \cite{bresar89}, and
\cite[Lemma 4.14]{CKLW} of Akemann and Pedersen \cite{akemann}, one does not need to assume $\cA$ or $\cB$ is unital.
Together with our new Theorem \ref{theorem:4.6}, which is a non-unital version of \cite[Theorem 4.6]{CKLW}, we can
now make use of the same proof of \cite[Theorem 4.11]{CKLW} to establish the assertion.
\end{proof}

Motivated by the theory of Banach lattices (see, e.g., \cite{abra}),
we call two C*-algebras being \emph{d-isomorphic} if there is a bounded bijective linear
map between them sending zero-products to zero-products.
We end this note with the following

\begin{cor}\label{cor:new2}
Two  C*-algebras are d-isomorphic if and only if they
are $*$-isomorphic.
\end{cor}

\begin{proof}
The conclusion follows from Theorem \ref{theorem:calgfull} and
a result of Sakai \cite[Theorem 4.1.20]{sakai} stating that
two algebraic isomorphic C*-algebras are indeed $*$-isomorphic.
\end{proof}

\end{document}